\documentclass[12pt]{article}
\usepackage{amstext}
\usepackage{graphicx}
\usepackage{caption}
\usepackage{subcaption}

\begin{document}
	\title{Analysis of a remarkable singularity in a nonlinear DDE }
	\author {Matthew Davidow$^\ast$, B. Shayak $^{\ast\ast}$ and {Richard H. Rand}$^{\ast\ast\ast}$}
	\date{}
		\maketitle 
{$^\ast$Center for Applied Mathematics, Cornell University, Ithaca, NY \\  
		$^{\ast\ast}$Sibley School of Mechanical and Aerospace Engineering, Cornell University, Ithaca, NY \\
		$^{\ast\ast\ast}$Department of Mathematics and Sibley School of Mechanical and Aerospace Engineering, Cornell University, Ithaca, NY }
	\section{Introduction}
In this work we investigate the dynamics of the nonlinear DDE (delay-differential equation)
\begin{equation}
\frac{d^2x}{dt^2}+x(t-T)+x^3 = 0
\label{1} 
\end{equation}
where $T$ is the delay.  Using Pontryagin's Principle, Bhatt and Hsu \cite{bhatt} showed that the origin in this equation is linearly unstable for 
all values of $T$ $>$ 0.  For $T=0$ however, the origin is obviously Liapunov stable.  Thus a stability change occurs as $T$ changes from zero to 
any  positive value, no matter how small. Associated with this change in stability is a remarkable bifurcation in which an infinite number of limit cycles exist for positive values of $T$ in the neighborhood of $T = 0$, their amplitudes going to infinity in the limit as $T$ approaches zero.\\

We investigate this situation in three ways:\\
1) Harmonic Balance,\\
2) Melnikov's integral,\\
3) Adding damping to regularize the singularity.\\

\section{Harmonic Balance}

We seek an approximate solution to eq.(\ref{1}) in the form:
\begin{equation}
x(t) = A \cos \omega t 
\label{2} 
\end{equation}

Substituting eq.(\ref{2}) in eq.(\ref{1}), simplifying the trig, and equating to zero the coefficients of $\sin \omega t$ and $\cos\omega t$ respectively, we obtain

\begin{equation}
 \sin \omega T = 0 \mbox{~~~~~~~~~~and ~~~~~~~~~~~}-\omega^2+\cos\omega T +\frac{3}{4} A^2 = 0
\label{3} 
\end{equation}

The first of these gives $\omega T=n\pi$ for $n$=1,2,3,$\cdots$, whereupon the second gives
\begin{equation}
A=\frac{2}{\sqrt{3}}\sqrt{\frac{n^2\pi^2}{T^2}\pm 1~}~,~~~~~~~~\mbox{$n$=1,2,3,$\cdots$}
\label{4} 
\end{equation}
where the upper sign refers to $n$ odd, and the lower sign refers to $n$ even.
For example, when $T=0.3$, Table 1 gives values for amplitudes of limit cycles for given values of $n$, from eq.(\ref{4}).

\begin{table}	
		\caption{Limit cycle amplitudes $A$ for values of $n$ in eq.(\ref{4}), for $T=0.3$.}
\hspace*{2.5 in}	
\begin{tabular}[t]{|c|c|}
		\hline
			$n$ & $A$\\
			\hline
		1& 12.14\\
		2& 24.21\\
		3& 36.29\\
		4& 48.38\\
		5& 60.47\\
		6& 75.56\\
		7& 84.65\\
		8& 96.74\\
		9& 108.83\\
		\hline
		\end{tabular}
\end{table}
Numerical integration of eq.(\ref{1}) using DDE23 in MATLAB shows limit cycles with amplitudes 12.31 and 33.56, which correspond to the approximate values 12.14 and 36.29 in Table 1.  See Fig.1.  Presumably the reason we do not see limit cycles with the other amplitudes listed in Table 1 is that they are unstable.  In fact, initial condition (x,x')=(26.681,0) for $t\leq 0$ leads to periodic motion with amplitude 12.14, while  initial condition (x,x')=(26.682,0) for $t\leq 0$ leads to periodic motion with amplitude 36.29, leading to the conclusion that there is an unstable periodic motion with amplitude approximately equal to 26.68, presumably corresponding to amplitude value 24.21 in Table 1.

\begin{figure}[h!]
	\begin{center}
		\includegraphics[width=0.9\textwidth]{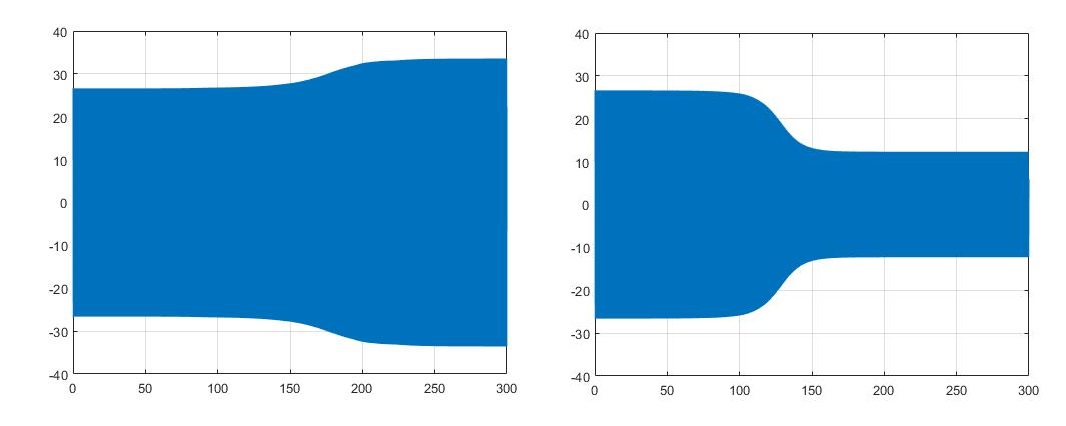}
		\caption{Numerical integration of eq.(\ref{1}) using DDE23 in MATLAB shows that
			initial condition (x,x')=(26.682,0) for $t\leq 0$ leads to periodic motion with amplitude 36.29 while  initial condition (x,x')=(26.681,0) for $t\leq 0$ leads to periodic motion with amplitude 12.14, leading to the conclusion that there is an unstable periodic motion with amplitude approximately equal to 26.68, presumably corresponding to amplitude value 24.21 in Table 1.}
		\label{fig:fig1}
	\end{center}
\end{figure}

\section{ Melnikov's integral}

We begin by generalizing the discussion to a wider class of systems, returning to eq.(\ref{1})
later.  We consider a conservative (Hamiltonian) system of the form:
\begin{equation}
\frac{dx}{dt}=\frac{\partial H}{\partial y},\mbox{~~~~~~}\frac{dy}{dt}=-\frac{\partial H}{\partial x}
\label{hom7}
\end{equation}
Note that eq.(\ref{hom7}) possesses the first integral $H(x,y)=$ constant, since $dH/dt=H_x \dot x+H_y\dot y=0$.\\

Now we add a perturbation to the conservative system (\ref{hom7}):
\begin{equation}
\frac{dx}{dt}=\frac{\partial H}{\partial y}+ g_1,\mbox{~~~~~~}\frac{dy}{dt}=-\frac{\partial H}{\partial x}
+ g_2
\label{hom8}
\end{equation}
where $g_1$ and $g_2$ are given functions of $x$ and $y$.\\

For the system (\ref{hom8}), the condition for one of the closed curves $H(x,y)=$ constant to be preserved under the perturbation (\ref{hom8}) turns out to be given by the vanishing of the following Melnikov integral:
\begin{equation}
\oint_\Gamma \left(g_1 \dot{y}-g_2 \dot{x}\right) dt=0
\label{hom9}
\end{equation}
where $\Gamma$ represents the unperturbed closed curve $H(x,y)=$ constant and where $\dot{x}$ and $\dot{y}$ refer to time histories around $\Gamma$ in the unperturbed system. The derivation uses Green's Theorem of the Plane, and the result is approximate (see section 3.3 in \cite{rhr}).\\

To apply the foregoing setup to eq.(\ref{1}), we write (\ref{1}) in the following form:
\begin{equation}\label{qq} \dot x=y \end{equation}
\begin{equation}\label{qq} \dot y=-x-x^3+(x-x(t-T)) \end{equation}
where $x$ written without an argument stands for $x(t)$.  That is we consider eq.(\ref{1}) to be a perturbed Hamiltonian system (\ref{hom8}) with
Hamiltonian
\begin{equation}
H(x,y)=\frac{1}{2}y^2+\frac{1}{2}x^2 +\frac{1}{4}x^4
\label{hom10}
\end{equation}
and with perturbations
\begin{equation}
 g_1=0 ~~~~\mbox{and}~~~~~ g_2=x-x(t-T)
\label{hom11}
\end{equation}

Thus in our case the Melnikov integral condition (\ref{hom9}) becomes
\begin{equation}
\int_0^P -(x(t)-x(t-T))\dot x(t) dt=\int_0^P x(t-T))\dot x(t) dt = 0
\label{hom12}
\end{equation}
where P is the period of the motion around $\Gamma$ in the unperturbed system, and where we have used the fact that:
 $$\int_0^P -x\dot x(t) dt=\left.\frac{x(t)^2}{2}\right|_0^P=\frac{x(P)^2-x(0)^2}{2}=0$$  
 Note that $x(P)=x(0)$ because $x(t)$ is periodic with period $P$.
Here $x(t)$ is the solution to eqs.(\ref{hom7}) with Hamiltonian (\ref{hom10})
 which turns out to be a Jacobian elliptic cn function, which may be written as  
 \begin{equation}
x=a_1\text{cn}(a_2t,k),
 \label{hom13}
 \end{equation}
 where the
 parameters $a_1$,$a_2$ and $k$ are related as follows (see section 2.2 in \cite{rhr}):
\begin{equation}
a_2^2=a_1^2+1,~~~~~ k^2=\frac{a_1^2}{2(1+a_1^2)}.
\label{hom14}
\end{equation}

Thus our Melnikov integral condition (\ref{hom12}) simplifies to:
\begin{equation}
\label{qqI} 
\int_0^P \text{cn}(a_2(t-T),k)~\frac{d}{dt}(\text{cn}(a_2t,k)~ dt=\int_0^P \text{cn}(a_2(t-T),k)~\text{sn}(a_2t,k)~\text{dn}(a_2t,k)~dt = 0
\end{equation}
where $P=4K(k)/a_2$, where $K(k)$ is a complete elliptic integral of the first kind.\\

In order to obtain an analytical approximation for this integral, we use the following expansions for the elliptic functions sn, cn and dn \cite{byrd}:\\
\begin{eqnarray}
\label{foo1}
\text{sn}(z,k)&=&\frac{2\pi}{kK}\sum_{n=0}^\infty\frac{q^{n+1/2}\text{sin}((2n+1)G)}{1-q^{2n+1}}\\
\label{foo2}
\text{cn}(z,k)&=&\frac{2\pi}{kK}\sum_{n=0}^\infty\frac{q^{n+1/2}\text{sin}((2n+1)G)}{1+q^{2n+1}}\\
\label{foo3}
\text{dn}(z,k)&=&\pi/(2K)+\frac{2\pi}{K}\sum_{n=0}^\infty\frac{q^{n}\text{cos}(2nG)}{1+q^{2n1}}
\end{eqnarray}
where  
$G=\pi z/(2K(k))$, $q=e^{-\pi K'(k)/K(k)}$ and $K'(k)=K(\sqrt{1-k^2})$.
We take the first term in each of the expansions (\ref{foo1}),(\ref{foo2}),(\ref{foo3}), whereupon the
Melnikov integral condition (\ref{qqI}) becomes:
\begin{equation}
\label{qqII} 
\int_0^P \cos(\pi a_2(t-T)/(2K))~\sin(\pi a_2t)/(2K))~dt = 0
\end{equation}
Expanding the cosine term gives
\begin{eqnarray}
\nonumber
\int_0^P [\text{sin}^2(\pi a_2 t/(2K))\text{sin}(\pi a_2 T/(2K))~~~~~~~~~~~~~~~~~~~~~~~~~~~~~~~~~~~~~~~~~~~~~~~\\
\label{qqIII} 
 +\text{sin}(\pi a_2 t/(2K))\text{cos}(\pi a_2 t/(2K))\text{cos}(\pi a_2T/(2K))] dt=0 
 \end{eqnarray} 

We are integrating over one full period, and thus the second term will integrate to 0. The first term, $\text{sin}^2(\pi a_2 t/(2K))$, is always positive and thus integrates to 0 only if the coefficient $\text{sin}(\pi a_2 T/(2K))$ is 0, i.e. eq.(\ref{qqIII}) becomes:
\begin{equation}
\label{qqIV} 
\text{sin}(\pi a_2 T/(2K))=0
\end{equation}

We are interested in the relationship between the amplitude $a_1$ and the delay $T$. The above gives an implicit relationship between $a_1$ and $T$ since $a_2^2=1+a_1^2$ and $K$ is  also determined by $a_1$ (through an elliptic integral). To make a much simpler explicit relationship we will use the fact that we are in the regime of $T<<1$, and in this parameter range we have empirically found that $a_1>>1$. Then from eqs.(\ref{hom14}) we can approximate $a_2 \approx a_1$, $k^2=a_1^2/(2a_2^2)\approx 1/2$. 

This gives us 

\begin{equation}
\label{qqV} 
\text{sin}(\pi a_1 T/(2K(1/2))=0 ~~\Rightarrow~~ a_1=2Kn/T
\end{equation}
where $K=K(1/2) \approx 1.854$, giving the result:

\begin{equation}
\label{qqVI} 
a_1 \approx 3.71~ n/T.
\end{equation}

This result may be compared to the Harmonic Balance result of eq.(\ref{4}), which is

\begin{equation}
\label{qqVII} 
a_1 \approx (2\pi/\sqrt{3}) n/T \approx 3.63 ~n/T.
\end{equation}

\begin{figure}
	\centering
	\includegraphics[width=0.7\textwidth]{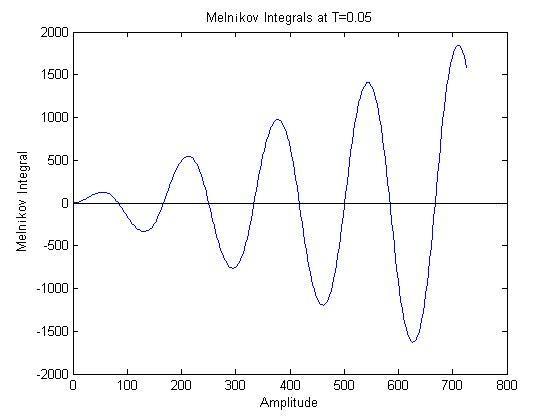}
	\caption{\label{fig:M5} Melnikov Integrals at T=0.05}
\end{figure}

\begin{figure}
	\centering
	\includegraphics[width=0.7\textwidth]{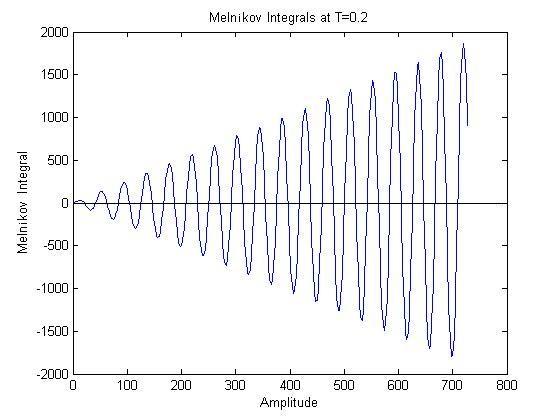}
	\caption{\label{fig:M2} Melnikov Integrals at T=0.2}
\end{figure}

\begin{figure}
	\centering
	\includegraphics[width=0.7\textwidth]{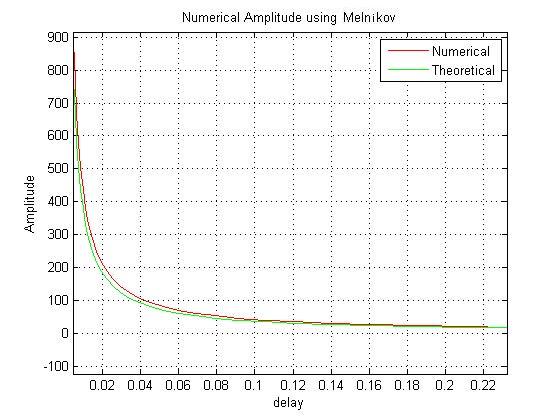}
	\caption{\label{fig:MA} Comparison of limit cycle amplitudes obtained numerically versus analytically.  Numerical values correspond to the first zero of the Melnikov integral (\ref{hom12}), while analytical values are those obtained by Harmonic Balance, eq.(\ref{qqVII}), for $n$=1.}
\end{figure}

These approximate analytical results may be checked by evaluating the Melnikov integral (\ref{qqI}) numerically.
For a fixed value of delay $T$, a value for the second integral in (\ref{qqI}) may be computed in MATLAB once the amplitude $a_1$ is chosen.
By varying $a_1$ we obtained two plots, one with delay $T=0.05$, and the other with $T=0.2$, see Figs.\ref{fig:M5} and \ref{fig:M2}. 
If we look at the zeros of both plots, it looks like they occur at integer multiplies of a certain amplitude. This agrees with the Harmonic Balance result 
of eq.(\ref{4}).  Fig.\ref{fig:MA} compares the numerical results with those of Harmonic Balance in a plot of the first zero (corresponding to $n=1$) for different values of delay. 

\section{Adding damping to regularize the singularity}

\begin{figure}[h!]
	\begin{center}
		\includegraphics[width=0.7\textwidth]{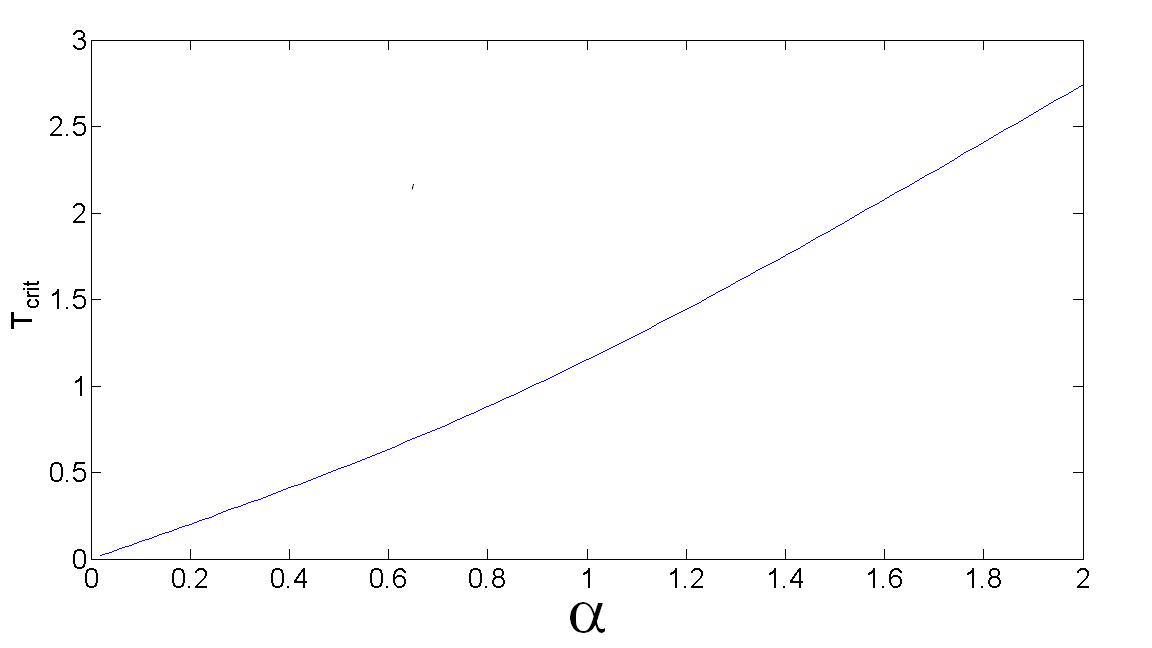}
		\caption{A plot of $T_{crit}$ from (\ref{d4})  as a function of $\alpha$. Note that at small $\alpha$ the function is like the identity $T=\alpha$. Limit cycles exist above this line, but not below it.}
		\label{fig:alpha1}
	\end{center}
\end{figure}

We have seen that in the case of eq.(\ref{1}), even infinitesimal delay gives rise to effective negative damping and growing oscillations. Accordingly, we expect that 
if damping is added to eq.(\ref{1}), as in the case of the following DDE:
\begin{equation}
\frac{d^2x}{dt^2}+\alpha \frac{dx}{dt}+x(t-T)+x^3 = 0
\label{d1} 
\end{equation}
then if $\alpha$ is held fixed and delay $T$ is increased from $0$, there will be a point at which the equilibrium at the origin will make a transition from stable to unstable.  Supposing that such a transition is a Hopf bifurcation, we linearize eq.(\ref{d1}) by dropping the $x^3$ term, and then set $x=\exp{i\omega t}$, giving the real and imaginary parts:
\begin{eqnarray}
\label{d2}
-{{\omega }^{2}}+\cos \omega T =0\\
\label{d3}
\alpha \omega -\sin \omega T =0
\end{eqnarray}   
Squaring and adding (\ref{d2}) and (\ref{d3}) and using (\ref{d2}) again yields the critical delay for a Hopf:
\begin{equation}
{{T}_{crit}}=\sqrt{2}~\displaystyle{\frac{\arccos 
		\displaystyle{
	\frac{-{{\alpha }^{2}}+\sqrt{{{\alpha }^{4}}+4}}{2}}}
{\displaystyle{{\sqrt{-{{\alpha }^{2}}+\sqrt{{{\alpha }^{4}}+4}}}}}}
\label{d4}
\end{equation}
A plot of $T_{crit}$ as a function of $\alpha$ can be seen in Fig.\ref{fig:alpha1}.

In addition to this Hopf bifurcation, it turns out that additional limit cycles can occur in this system by being born in a fold (also known as a saddle-node of cycles).  In order to see this we again use the method of Harmonic Balance.
Assuming an approximate solution of the form $x(t) = A \cos \omega t$, we substitute into eq.(\ref{d1}), simplify the trig, and equate to zero the coefficients of $\sin \omega t$ and $\cos\omega t$ respectively, giving:

\begin{equation}
\sin \omega T = \alpha \omega \mbox{~~~~~~~~~~and ~~~~~~~~~~~}-\omega^2+\cos\omega T +\frac{3}{4} A^2 = 0
\label{d5} 
\end{equation}

\begin{figure}[h!]
	\begin{center}
		\includegraphics[width=0.7\textwidth]{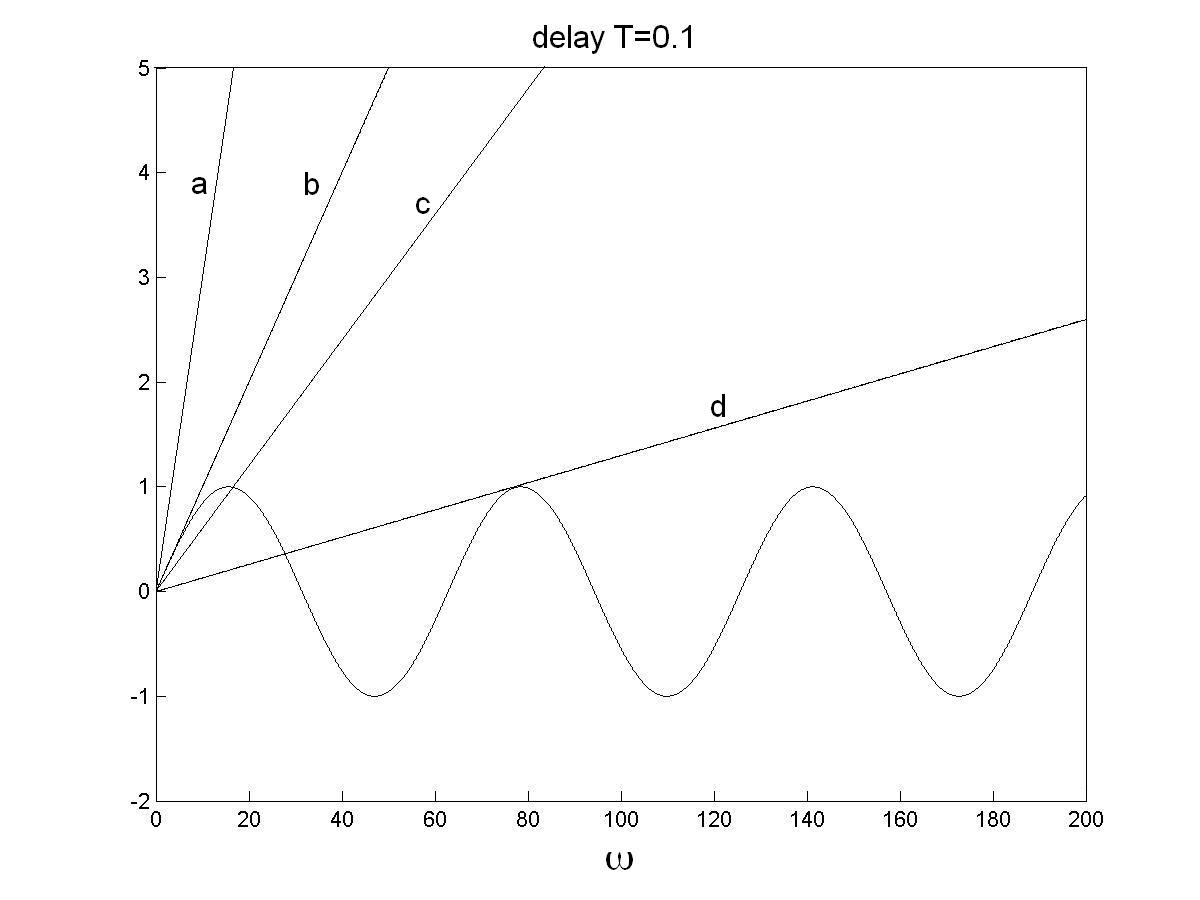}
		\caption{Graphical representation of the first of eqs.(\ref{d5}).  The straight lines are $y=\alpha\omega$ and have slope $\alpha$. 
			The sinusoid is $y=\sin \omega T$.}
		\label{fig:foo1}
	\end{center}
\end{figure}

Suppose the value of $T$ is fixed and $\alpha$ is started from a high value. The first equation of (\ref{d5}) can be viewed in terms of two functions of the variable $\omega$; one is the straight line $\alpha \omega$ and the second is the sinusoid $\sin \omega T$. 
See Fig.\ref{fig:foo1}.\\

If $\alpha > T$ then the two curves intersect only at the trivial point $\omega=0$ and there is no limit cycle. This corresponds to curve $a$ in Fig.\ref{fig:foo1}.  Now consider the situation as $\alpha$ is lowered. When it becomes equal to $T$ there is a tangency at the origin (curve $b$ in Fig.\ref{fig:foo1}), and upon being slightly lower still, the curves develop a non-trivial intersection (curve $c$ in Fig.\ref{fig:foo1}). This means that a sinusoidal response with frequency $\omega$ is a possible state of the system. Once the frequency is specified, the amplitude of the motions gets determined by the second equation in (\ref{d5}). We thus have a limit cycle with amplitude $A$ and frequency $\omega$. \\

As we lower $\alpha$ still further, the non-trivial intersection point between the straight line and the sinusoid will shift rightwards. Ultimately, the two graphs will touch at a second point (curve $d$ in Fig.\ref{fig:foo1}) and a new pair of limit cycles will get born there, since further lowering $\alpha$ will turn the tangency into a pair of intersections, one corresponding to a stable limit cycle and the other to an unstable one. Thus, we can say that a saddle node bifurcation of cycles is occurring there. The points of tangency are given by the relation
\begin{equation}
\sin \omega T =\alpha \omega
\end{equation}
\begin{equation}
T \cos \omega T =\alpha
\end{equation}
which imply that at the $n^{th}$ intersection point
\begin{equation}
\omega =\frac{1}{T }{{\beta }_{n}}
\label{goo1}
\end{equation}
\begin{equation}
{{\alpha }_{n}}=T \cos {{\beta }_{n}}
\label{goo2}
\end{equation}
where $\beta _{n}$'s are the solutions of $\tan x = x$.\\

\begin{figure}[h!]
	\begin{center}
		\includegraphics[width=0.7\textwidth]{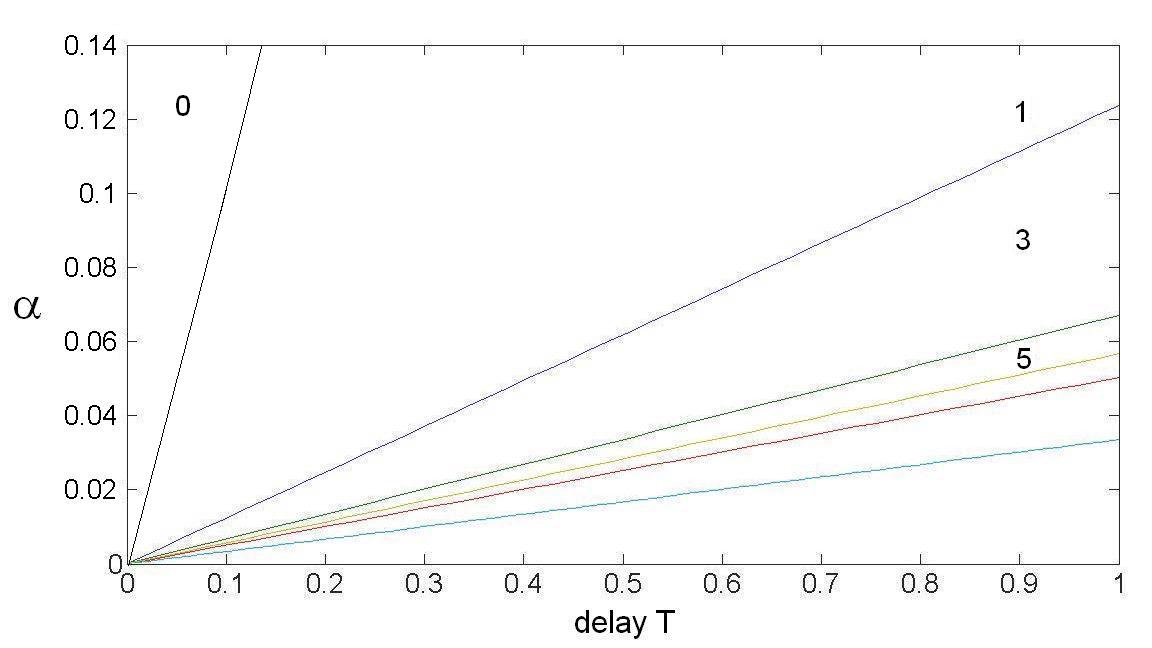}
		\caption{Plot of the bifurcation curves using the tangency condition (\ref{goo1}),(\ref{goo2}). The number of limit cycles in the various regions is shown in the first few cases (0,1,3,5). }
		\label{fig:foo2}
	\end{center}
\end{figure}

A bifurcation diagram using these relations is shown in Fig.\ref{fig:foo2}.

In this figure, the line on the left is the Hopf bifurcation. The uppermost line on the right is the first saddle-node bifurcation of cycles, while the next lines going down correspond to the subsequent saddle-node bifurcations. These predictions are in agreement with numerical simulation results. A comparison between theory and simulation 
 performed on Matlab using the routine DDE23,
 is presented in Table 1.  This Table shows the birth of a limit cycle (LC) as $\alpha$ is lowered to a value smaller than $T$ and its branching out into multiple cycles as $\alpha$ is lowered further.\\

\begin{table}
\begin{tabular}{c|c|c|c|c}
	\hline
	$T$ & $\alpha$ & Eigenvalues & Calculated LC Amplitude & Observed LC Amplitude\\
	\hline
	0.4 & 0.1 & $0.13 \pm  0.95i$ & 7.29 & 7.3 \\
	& 0.2 & $0.09 \pm  0.90i$ & 5.64 & 5.5 \\
	& 0.3 & $0.04 \pm  0.96i$ & 3.70 & 3.6 \\
	& 0.4 & NRP & DNE & DNE \\
	\hline
	0.6 & 0.1 & $0.20 \pm  0.91i$ & 5.25 & 5.4 \\
	& 0.2 & $0.16 \pm  0.95i$ & 4.54 & 4.6 \\
	& 0.3 & $0.12 \pm  0.95i$ & 3.76 & 3.7 \\
	& 0.4 & $0.07 \pm  0.95i$ & 2.80 & 2.8 \\
	& 0.5 & $0.02 \pm  0.95i$ & 1.78 & 1.8 \\
	& 0.6 & NRP & DNE & DNE \\
	\hline
	1.0 & 0.1 & & 3.59, 8.57, 9.91 & 3.6, 10 \\
	2.0 & 0.1 & & 2.07, 3.64, 5.40, 7.40, 8.59 & 2.1, 5.4, 8.8\\
	\hline
\end{tabular} 
\caption{ 
	This Table shows the theoretically calculated and numerically observed amplitudes of the LCs as parameter values are varied. (NRP=negative real part, DNE=does not exist.) In some regions of the space, 3 and 5 LCs are seen in the harmonic balance. In this case, the first one and then the subsequent alternate ones are observed numerically with the intermediate amplitudes acting as separatrix. It is also seen that the LC is born when the real part of the eigenvalue crosses from negative to positive.}
\end{table}

\newpage
\section{Conclusions}

We have investigated the occurrence of limit cycles in the delay-differential equation:

\begin{equation}
\frac{d^2x}{dt^2}+x(t-T)+x^3 = 0
\label{1a} 
\end{equation}

Besides numerical integration, we used three different approximate analytic approaches to study this system.
All of these approaches support the conclusion that this system exhibits an infinite number of limit cycles for positive values of $T$ in the neighborhood of $T = 0$, their amplitudes going to infinity in the limit as $T$ approaches zero.

\end{document}